\documentclass[twocolumn]{webofc}

\def\RR{\mathbb{R}}

\def\1{\mathbf{1}}
\def\<{\langle} \def\>{\rangle}

\usepackage{array} %% needed for advanced table manipulation
\newcolumntype{L}[1]{>{\raggedright\let\newline\\\arraybackslash\hspace{0pt}}m{#1}}
\newcolumntype{C}[1]{>{\centering\let\newline\\\arraybackslash\hspace{0pt}}m{#1}}
\newcolumntype{R}[1]{>{\raggedleft\let\newline\\\arraybackslash\hspace{0pt}}m{#1}}

\usepackage[varg]{txfonts}   % Web of Conferences font
\usepackage{amsmath,amssymb}
\usepackage{graphicx}
\usepackage{algorithm}
\usepackage{algpseudocode} 
\usepackage[usenames,dvipsnames]{color}
\usepackage{colortbl,tikz}
\usepackage{url}

\begin{document}
\title{A fuzzy-set theoretical framework for computing exit rates of rare events in potential-driven diffusion processes}
%
% subtitle is optionnal
%
%%%\subtitle{Do you have a subtitle?\\ If so, write it here}

\author{\firstname{Marcus} \lastname{Weber}\inst{1}\fnsep\thanks{\email{weber@zib.de}} \and
        \firstname{Natalia} \lastname{Ernst}\inst{1,2}
      }

\institute{Zuse Institute Berlin (ZIB), Takustr. 7, 14195 Berlin, Germany
\and
           Technical University Berlin (TUB), Straße des 17. Juni 135, 10623 Berlin, Germany
          }

\abstract{%
  This article is about molecular simulation. However, the theoretical results apply for general overdamped Langevin dynamics simulations. Molecular simulation is often used for determining the stability of a complex (e.g., ligand-receptor). The stability can be measured by computing the expected holding time of the complex before its dissociation. This dissociation can be seen as an exit event from a certain part $S$ of the conformational state space $\Gamma$. Determining exit rates (i.e, for SDE-based simulations exiting from a given starting set $S$) for a stochastic process in which the exit event occures very rarely is obviously hard to solve by straight forward simulation methods. Finding a low variance procedure for computing rare event statistics is still an open problem. Imagine now, e.g., a simulation of a diffusion process. As long as the time-dependent state trajectory is inside the starting set $S$, no information is gained about the rare event statistics. Only at that point of time, when the process leaves the starting set, a piece of information about the exit rate is collected. If $S$, however, is a fuzzy set given by a membership function, then there might be additional information of the kind ``the process is about to leave the set''. However, how to define an exit rate from a fuzzy set? 
}
\maketitle
%\doublespacing
%

\section{Introduction}
\label{sec:intro}

\begin{figure}[ht]
\centering
\includegraphics[width=0.45\textwidth]{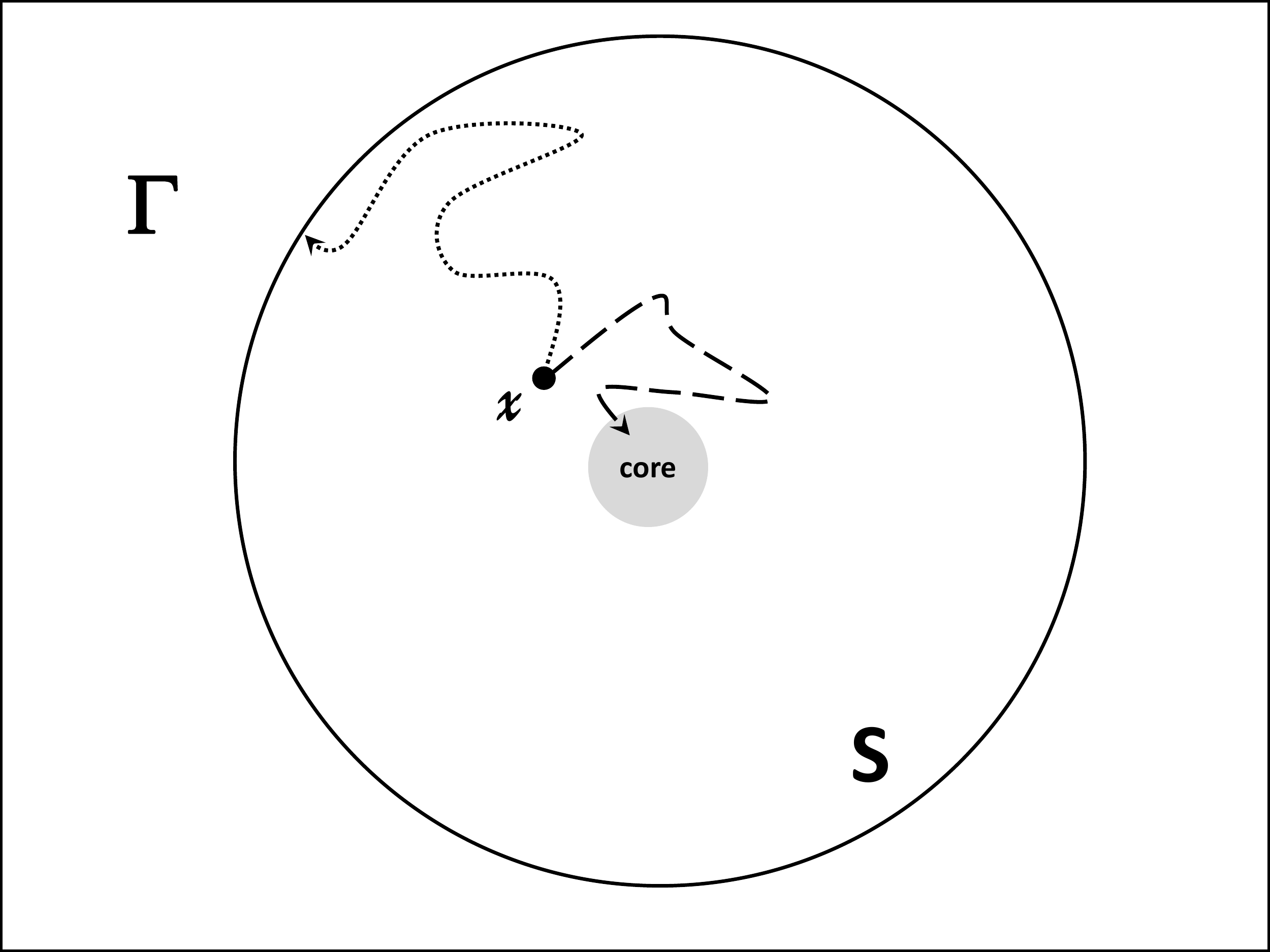}
\caption{\label{fig:situation} In the state space $\Gamma$ (rectangle) there is a subset $S$ of starting points (big circle). If we start trajectories of a stochastic process at a starting point $x\in S$ of a fixed time-length $t$, then there is a certain probability that we get a trajectory that is always within $S$ (dashed line) and that does not leave $S$ (dotted line). We are interested in how this probability depends on the time-length $t$ and on the starting point $x$. The ``core'' is a subset of $S$ which is very likely to be reached by trajectories starting in $S$.}
\end{figure}
In molecular simulation we are often faced with the situation that is depicted in Fig.~\ref{fig:situation}. Let us assume, that we run a Brownian dynamics (potential-driven diffusion, overdamped Langevin) simulation guided by a given potential energy. This is mathematically formulated in Eq.(\ref{diffusion}) below. Brownian dynamics trajectories that start in low-energy regions $S$  of the potential energy surface will leave these regions (dotted line) only with a very small probability. It is more likely to observe that the trajectory will approximate (dashed line) and dominantly sample states within the lowest energy part of $S$. This part is denoted as ``core'' in Fig.~\ref{fig:situation}. That means that it is very probable to observe trajectories to approximate the core set rather than trajectories which leave the low-energy region $S$. This ``exit event'' is, thus, a rare event. The (biological or chemical) stability of many molecular systems is given by the rareness of these events. A straight forward molecular simulation approach to estimate the exit rate is the following: We start several trajectories from $x\in S$ and determine the time they need to exit from $S$. This approach will be extremely inefficient, because the statistics will depend on rare events. 

Many approaches in literature try to overcome the sampling problem of rare events. Here are some examples:
\begin{itemize}
\item An uncountable number of molecular simulation methods try to accelerate the transitions between the molecular conformations. Girsanov's Theorem can be applied in order to reweight the accelerated samplings towards the original slow process \cite{SIAM14}. However, evaluating the reweighting formula and finding proper acceleration methods can be difficult. 
\item Computing rare events by milestoning \cite{Voronoi} discretizes the exit path of the molecular system and determines the (fast) rates between the milestones. However, this approach is based on a discretization scheme for high-dimensional spaces.  
\item Parallel Trajectory Splicing \cite{parallel} is a method that extensively makes use of parallel computation facilities on supercomputers. Since this method is not based on high-dimensional discretization schemes it is a promising way, if supercomputers are available to run many trajectories in parallel. Parallel Trajectory Splicing, however, is depending on (set-based) defining different states of the system. It also depends on a certain eigenvalue gap between two consecutive eigenvalues of the infinitesinaml generator $\cal L^\ast$. This gap controls the error of the artificially generated long-time simulation (from many short-time trajectories). 
\end{itemize}

What all approaches have in common, is that they try to estimate an exit rate out of a given (metastable) set $S$. Reaching the boundary of that set $S$ during a Brownian dynamics simulation is a rare event. If $S$ is a set, then there is in principle no preferences between ``approaching the core'' (Fig.~\ref{fig:situation}, dashed line) or ``being about to leave the set $S$'' in these approaches. There is no further discrimmination between those trajectories which stay inside $S$. There is only a discrimmination between trajectories which either leave the set or stay in $S$. 

This article will start with defining a membership function (fuzzy set) $\chi$ instead of $S$ by applying the Robust Perron Cluster Analysis~ \cite{PCCAplus, weber} (PCCA+) to the eigenfunctions of some infinitesimal generator. This fuzzyness will provide the discussed discrimmination.

%%%%%%%%%%%%%%%%%%%%%%%%%%%%
\section{Theoretical Background}\label{sec:theo}

\subsection{The stochastic differential equation}\label{sec:process}
The presented theory is based on reversible, ergodic, and stochastic equations of motion \cite{CLN}. The motion is a diffusion in an energy landscape $V:\Gamma \rightarrow\RR$, i.e., the realizations of the dynamics lead to trajectories $x_t\in\Gamma, t\in[0,\infty),$ with 
\begin{equation}\label{diffusion}
d\, x_t =-\nabla_x V(x_t)\,dt + \sigma \, dB_t,
\end{equation}
where $B_t$ denotes Brownian motion and $\sigma$ the constant diffusion parameter. The non-linear diffusion equation (\ref{diffusion}) generates trajectories of states. If instead of a single state $x_t$ a density of states $\pi_t$ is propagated with the above dynamics, then the corresponding linear Fokker-Planck equation describes this propagation of densities \cite{handbook}: 
\begin{equation}\label{fokkerplanck}
\frac{\partial\pi_t}{\partial t} = {\cal L}\pi_t = -\mathrm{div}_x(\nabla_x V(x)\,\pi_t) + \frac{\sigma^2}{2}\Delta_x \pi_t,
\end{equation}
where $\cal L$ is the corresponding Fokker-Planck operator and where $\Delta_x$ denotes the Laplacian operator. We assume ergodicity of this process, which leads to the invariant density $\pi$ defined by ${\cal L}\pi=0$. $\pi$ gives rise to a weighted scalar product $\langle \cdot, \cdot \rangle_\pi$. Note that, although the diffusion equation (\ref{diffusion}) is non-linear, the corresponding propagation of densities is given by a linear equation (\ref{fokkerplanck}). In the following, the adjoint operator ${\cal L}^+$ of $\cal L$, i.e., the infinitesimal generator of the stochastic process in the form
\[
{\cal L}^+=-\nabla_x V(x)\cdot \nabla_x + \frac{\sigma^2}{2}\Delta_x,
\]
will play an important role. We will use ${\cal L}^\ast=-{\cal L}^+$. From this infinitesimal generator the transfer operator ${\cal P}^\tau$ of the process for the time-length
$\tau$ can be derived as
\begin{equation}\label{eq:generator}
 {\cal P}^\tau=\exp(-\tau{\cal L}^\ast).
\end{equation}
The transfer operator is a widely used tool to compute transition probabilities between molecular conformations \cite{chhabil}. Its spectral properties \cite{A19-6} are strongly connected to the spectral properties of $\cal L$ and ${\cal L}^\ast$. 

\subsection{Defining membership functions}\label{sec:PCCA+}
The algorithm PCCA+~ \cite{PCCAplus, weber} is commonly used for identifying metastable membership functions $\chi$ of molecular systems. The way this algorithm works is the following: Eigenfunctions of the transfer operator ${\cal P}^\tau$ are computed which correspond to eigenvalues close to the eigenvalue $\lambda=1$ of ${\cal P}^\tau$.   PCCA+ defines membership functions as a linear combination of these eigenfunctions. Due to (\ref{eq:generator}), the eigenfunctions of ${\cal P}^\tau$ correspond to the eigenfunctions of ${\cal L}^\ast$. The constant function ${1\!\! 1}$ is always the first (dominant) eigenfunction. In the easiest case, $\chi$ is a non-trivial linear combination of ${1\!\! 1}$ and an eigenfunction $f:\Gamma\rightarrow \RR$ of ${\cal L}^\ast$. Let $f$ be an eigenfunction of ${\cal L}^\ast$ corresponding to some eigenvalue $\bar{\epsilon}\not=0$. 
From the reversibility of the process one can derive that $\{{1\!\! 1}, f\}$ are orthonormal functions (with regard to the scalar product $\langle\cdot,\cdot\rangle_\pi$). Then $\chi$ is a linear combination of these two functions with $\chi=\bar{\alpha}f+\bar{\beta}{1\!\! 1}$ with real numbers $\bar{\alpha},\bar{\beta}$. This provides:
\[
f=\frac{1}{\bar{\alpha}}\chi-\frac{\bar{\beta}}{\bar{\alpha}}{1\!\! 1}.
\]
Thus,
\begin{equation}\label{eq:shiftscale}
{\cal L}^\ast \chi=\bar{\alpha}\bar{\epsilon} f = \bar{\epsilon}\chi -\bar{\epsilon}\bar{\beta} {1\!\! 1}=\alpha \chi + \beta,
\end{equation}
where $\alpha=\bar{\epsilon}$ and $\beta=-\bar{\epsilon}\bar{\beta}$.

\subsection{Implication of PCCA+ result}\label{sec:comp}
Equation (\ref{eq:shiftscale}) can be understood in the following sence. The action of the infinitesimal generator ${\cal L}^\ast$ in the situation of Sec.~\ref{sec:PCCA+} is like shifting and scaling the function $\chi$. If equation (\ref{eq:shiftscale}) holds, what does it imply? 
Equation (\ref{eq:shiftscale}) is equivalent to
\begin{equation}\label{combine2}
{\cal L}^\ast \chi = \epsilon_1 \chi - \epsilon_2 (1-\chi),  
\end{equation}
where $\epsilon_1=\alpha + \beta$ and $\epsilon_2=-\beta$ are fixed  numbers defined by the PCCA+ algorithm. By multiplying (\ref{combine2}) with the expression $e^{-\epsilon_1 t}$, and by defining a function $p_\chi(x,t) := \chi(x) \,e^{-\epsilon_1 t}$, we equivalently get:
\begin{equation}\label{combine}
-{\cal L}^\ast p_\chi - \epsilon_2 \,\frac{1-\chi}{\chi}\,p_\chi = -\epsilon_1 p_\chi,
\end{equation}
for all $x\in\Gamma$ with $\chi(x)\not=0$. The definition of $p_\chi$ can be expressed by an ordinary differential equation:
\begin{equation}\label{ratedef2}
\frac{\partial p_\chi}{\partial t}= -\epsilon_1 p_\chi \quad \mathrm{and} \quad p_\chi(x,0)=\chi(x). 
\end{equation} 
Combining (\ref{ratedef2}) and (\ref{combine}) leads to the following differential equation:
\begin{equation}\label{holdingsoft2}
\frac{\partial p_\chi}{\partial t}=-{\cal L}^\ast p_\chi - \epsilon_2 \,\frac{1-\chi}{\chi}\,p_\chi, \quad p_\chi(x,0)=\chi(x).
\end{equation}
The equation (\ref{holdingsoft2}) is the solution of the following conditional expectation value problem according to the Feynman-Kac formula (Equations III.1 and III.2 in \cite{gzyl}):
\begin{equation}\label{holdingsoft}
p_\chi(x,t) = \mathbb{E}\left[\chi(x_t)\cdot \exp\Big(-\epsilon_2\int_0^t\frac{1-\chi(x_r)}{\chi(x_r)}\,dr\Big)\right]_{x_0=x}.
\end{equation}
In this expression, $x_t\in \Gamma$ are realizations of the stochastic differential equation starting in $x_0=x$.

\subsection{Holding probability of membership functions}
\label{sec:relax}
If we assume that $\exp(-\infty)=0$, then the equation (\ref{holdingsoft}) turns into  
\begin{equation}\label{holdingprob}
p_{1\!\! 1_S}(x,t)=\mathbb{E}\left[{1\!\! 1}_S(x_t)\cdot \delta_0\Big(\int_0^t(1-{1\!\! 1}_S(x_r))\,dr\Big)\right]_{x_0=x}.
\end{equation} 
Here, we replaced $\chi$ with a characteristic function ${1\!\! 1}_S$ of a set $S$. The characteristic function ${1\!\! 1}_S:\Gamma\rightarrow\{0,1\}$ of the set $S$, is only $1$ for $x\in S$. In equation (\ref{holdingprob}), the function $\delta_0:\RR\rightarrow \{0,1\}$ is only $1$ if the integral is zero. Equation (\ref{holdingprob}) is the well-established way of defining holding probabilities. Each realization of the stochastic process provides a trajectory $x_t\in \Gamma$. The starting set $S \subset \Gamma$ is a connected, open subset in the state space $\Gamma$. The {\em holding probability} $p_{{1\!\! 1}_S}(x,t)$ of the set $S$ is the percentage of realizations of the stochastic process starting in $x_0=x$ which have never left the set $S$ until time $t$, i.e., $x_r\in S$ for all $r\in[0,t]$. The longer the time $t$ the more trajectories will leave $S$, thus, the smaller $p_{{1\!\! 1}_S}(x,t)$. Equation (\ref{holdingsoft}) can, therefore, be seen as the definition of a {\em $\chi$-holding probability} of a membership function $\chi$. The term $\chi$-holding probability is different from the set-based definition of a holding probability. In this article we will define further quantities of this fuzzy-set-based kind. Under the assumptions of Sec. \ref{sec:comp} the $\chi$-holding probability decreases exponentially with time
\begin{equation}\label{ratedef}
p_\chi(x,t) = \chi(x) \,e^{-\epsilon_1 t}.
\end{equation}
From the structure of equations (\ref{holdingsoft}) and (\ref{ratedef}), we can derive that $\epsilon_1>0$ is the {\em $\chi$-exit rate} of the membership function $\chi$. Note that, equation (\ref{ratedef}) is consistent with the definition of $\chi(x)=p_\chi(x,0)$, because $\chi(x)$ is the probabilty to classify $x$ as a starting point of the diffusion process. Usually, one can not show that the set-based holding probability depends exponentially on the simulation time $t$. The following would just be an approximation:
\begin{equation}\label{approx}
   p_{{1\!\!1}_S}(x,t)\approx {1\!\! 1}_S(x)\, e^{-\epsilon t}.
\end{equation}   
Since (\ref{approx}) is not an equality, the definition of a set-based exit rate $\epsilon$ is often understood as a fitting parameter of exit time distributions \cite{srebnik2007}, as an optimization quantity \cite{hartmann2014}, or as an asymptotic value, e.g., in the large deviation principle \cite{cramer}. The membership-based $\chi$-holding probability decreases exponentially, if the conditions of Sec.~\ref{sec:comp} are satisfied.

\subsection{Formula for exit rates}\label{sec:formula}
According to the previously identified relations between the defined quantities $\epsilon_1$ and $p_\chi(x,t)$, we can define what an exit rate out of a fuzzy set should be.

{\bf Definition:} Given a membership function (fuzzy set) $\chi:\Gamma \rightarrow [0,1]$ for the starting points of a potential-based diffusion process with infinitesimal generator ${\cal L}^\ast$ such that ${\cal L}^\ast \chi = \alpha \chi + \beta {1\!\! 1}$, then $\epsilon_1=\alpha+\beta$ is the $\chi$-exit rate out of $\chi$, i.e., the $\chi$-holding probability meets $p_\chi(x,t) = \chi(x) \,e^{-\epsilon_1 t}$. 

\subsection{Formula for exit paths}\label{sec:paths}
If one wants to compute the exit path direction from a given state $x\in \Gamma$, then there are in principle two ways to define them in the $\chi$-context. First, the exit path is connected to a decreasing holding probability. Given the $\chi$-holding probability $p_\chi(x,t)$, the exit path direction can be defined as the negative gradient of $p_\chi$, because this is the direction starting in $x$ in which the holding probability decreases the most. Due to the formula of the exit rate, this {\em direction} is given by $-\nabla_x\chi(x)$. Second, given a membership function $\chi$, the exit path direction can be defined by a decreasing $\chi$-value. Both definitions lead to the same direction $-\nabla_x\chi(x)$ of the path. Following the negative gradient of $\chi$ from a given starting point $x$ provides the $\chi$-exit path. In the PCCA+ context, note that the gradient of $\chi$ and the gradient of the eigenfunction $f$ are linearly dependent, which means that also the eigenfunctions of ${\cal L}^\ast$ provide the exit paths in the situation of Sec.~\ref{sec:PCCA+}.   
 
\subsection{Mean holding time}
According to the theory \cite{Pav14} of stocahstic differential equations, the mean holding time of a process can be expressed by the integral of the holding probability. In our case, a corresponding definition of a {\em $\chi$-mean holding time} $ t_1(x)$ depending on the initial state $x_0=x$ of the process would be the following integral: 
\[
t_1(x)=\int_0^\infty p_\chi(x,t)\, dt=\chi(x)\,\frac{1}{\epsilon_1}. 
\]
Due to the choice of the sign of ${\cal L}^\ast$, one would further expect from theory \cite{Pav14} that a set-based mean holding time $t(x)$ meets ${\cal L}^\ast t(x)=1$ inside the open set $S$. We replaced $S$ by the fuzzy set $\chi$. For the $\chi$-mean holding time we get:
\[
{\cal L}^\ast t_1(x) = \frac{1}{\epsilon_1}{\cal L}^\ast\chi(x) = \chi(x) - \frac{\epsilon_2}{\epsilon_1} (1-\chi(x)).
\]
This is indeed what we expect, if $x\in S$ corresponds to $\chi(x)\approx 1$.
 
\subsection{Summarizing the main equations}
In the center of discussions there is a membership function $\chi$ which satisfies the ``almost'' eigenvalue equation
\[
{\cal L}^\ast \chi = (\epsilon_1+\epsilon_2)\chi - \epsilon_2,
\] 
with $\epsilon_2\ll \epsilon_1$ and some infinitesimal generator ${\cal L}^\ast$. The exit paths out of those fuzzy sets $\chi$ are given by following the gradients $-\nabla \chi(x)$.
On the basis of $\chi$, a new quantity denoted as $\chi$-holding probability is defined as 
\[
p_\chi(x,t) = \mathbb{E}\left[\chi(x_t)\cdot \exp\Big(-\epsilon_2\int_0^t\frac{1-\chi(x_r)}{\chi(x_r)}\,dr\Big)\right]_{x_0=x}.
\] 
Under certain conditions this quantity depends exponentially on time, i.e., 
\[
p_\chi(x,t)=\chi(x)e^{-\epsilon_1 t}.
\]
The positive number $\epsilon_1>0$ is the $\chi$-exit rate, leading to the $\chi$-mean holding time $t_1(x)=\chi(x)/\epsilon_1$, which is result of a partial differential equation:
\[
{\cal L}^\ast t_1(x) = \chi(x) - \frac{\epsilon_2}{\epsilon_1}(1-\chi(x)).
\]
This equation is ``almost'' like the equation ${\cal L}^\ast t(x)=1$ for computing set-based mean holding times in the interior of a given set $S$.

%%%%%%%%%%%%%%%%%%%%%%%%%%%%
\section{Approximations}\label{sec:algo}

\subsection{Implications of $\chi(x)\approx 1$}\label{sec:meaning}
The definition of a $\chi$-holding probability (\ref{holdingsoft}) is different from the set-based definition (\ref{holdingprob}). The random variable in the set-based definition can only be $0$ or $1$. Thus, only the ``exit event'' defines the $\chi$-exit rate. In the case of the $\chi$-holding probability the function $\chi$ is not constantly $1$ inside the fuzzy set of starting points. The $\chi$-exit rate could origin from the exponential penalty term
\[ 
   \exp\Big(-\epsilon_2\int_0^t\frac{1-\chi(x_r)}{\chi(x_r)}\,dr\Big).
\] 
Thus, maybe not the ``exit event'' produces the $\chi$-exit rate, but it stems from the fuzzy definition of $\chi$. We have to answer the question of the problematic time-scale, i.e., at what time $t_2$ the exponential expression starts to dominate the function $\chi$? In order to compute this, we will assume that $\chi$ is constant. The question is now, when will
\[
   \exp\Big(-\epsilon_2 \,t_2\, \frac{1-\chi}{\chi}\Big)\leq \chi ?
\]
This is the case for
\[
   t_2 \geq \frac{-\ln(\chi)\chi}{(1-\chi)\epsilon_2}.
\]
If we assume, that $\chi$ is nearly $1$ and take the limit of that expression for $\chi\rightarrow 1$, then the result is $t_2\geq \epsilon_2^{-1}$ according to the rule of De L'Hospital. This means, that the mean holding time of $\chi$ should be smaller than $\epsilon_2^{-1}$, such that the ``exit event'' dominates the definition of $p_\chi$. The mean holding time $t_1$ (identical to the mean first exit time) is given by $t_1=\chi(x)\,\frac{1}{\epsilon_1}$. Again we assume $\chi(x)\approx 1$. Thus, if $\epsilon_1 > \epsilon_2$, then the ``exit event'' dominates the definition of the holding probability. In this case,  the definition of the $\chi$-exit rate is consistent with our physical interpretation.  A lemma (Lemma 3.6 in \cite{weber}) provides the following connection between $\bar{\beta}$ in equation (\ref{eq:shiftscale}) and the statistical weight $\pi_\chi$ of $\chi$:
\[
\bar{\beta}=\pi_\chi:=\frac{\int_\Gamma \chi(x)\, \pi(x)\, dx}{\int_\Gamma \pi(x) \, dx},
\]  
where $\pi$ is the invariant density of the stochastic process. 
Using the easy formula of Sec~\ref{sec:formula}, the $\chi$-exit rate is computed as  $\epsilon_1=\bar{\epsilon}(1-\bar{\beta})=\bar{\epsilon}(1-\pi_\chi)$. Furthermore, $\epsilon_2=\bar{\epsilon}\pi_\chi$. The condition $\epsilon_1> \epsilon_2$ means that the concept of the $\chi$-holding probability and the corresponding $\chi$-exit rate is only physically meaningful, if $\pi_\chi< 0.5$, i.e., it is meaningful for starting points out of a small ``subset'' $\chi$ of the state space $\Gamma$.

\subsection{Approximating $\chi$}\label{sec:complex}
Section \ref{sec:formula} provides an equation for computing $\chi$-exit rates, which are physically meaningful under certain conditions mentioned in the last section. However, if we apply Sec.~\ref{sec:formula} in the PCCA+ context, then computing the $\chi$-exit rate is given by $\epsilon_1=\bar{\epsilon}(1-\pi_\chi)$. In order to compute one (let us denote it as) local property $\epsilon_1$, one needs two global properties of the system, namely one eigenvalue $\bar{\epsilon}$ of ${\cal L}^\ast$ and the statistical weight of $\chi$. These two quantities are correlated via the eigenfunction $f$. In practise, the problem of computing a $\chi$-exit rate would turn into a function approximation problem for $f$ in high-dimensional spaces $\Gamma$. Function approximation has in general a non-polynomial complexity. A lot of effort has been spent in order to circumvent this ``curse of dimensionality''. Approximations of eigenfucntions of the transfer operator ${\cal P}^\tau$ have been computed using Markov State Modeling \cite{MSMbuilder, Chodera07, SarichNoeSchuette_MMS10_MSMerror}, diffusion maps \cite{diffusionmaps}, the variational principle \cite{A19-6}, committor functions \cite{milestoning}, and many other mathematical tools. Note that the core set approach \cite{committor} approximates eigenfunctions of ${\cal P}^\tau$ based only on committor values which can be sampled by generating an ensemble of trajectories according to (\ref{diffusion}). 

One result of Sarich \cite{SarichPhD} (Theorem 13) is the following: The subspace spanned by two eigenfunctions $\{{1\!\! 1}, f\}$ and the subspace spanned by a committor function $\xi$ and ${1\!\! 1}-\xi$ almost coincide using suitably defined core sets. This is especially the case, if the stochastic process between the core sets is metastable. Thus, instead of computing $f$, one can compute a committor function $\xi$. Since $\xi$ is a function between $0$ and $1$, the PCCA+ result on the basis of this approximation space $\{\xi,{1\!\! 1}-\xi\}$ would be $\chi=\xi$. Thus, the membership function defined in Sec.~\ref{sec:PCCA+} is almost identical to a committor function $\xi$. 

However, how can we estimate a committor function, if there is only one core set? We invert the argument of the introduction and note that if we generate a trajectory of certain length $T$ (much smaller than the mean exit time) starting in $x$ and reach the pre-defined core within or before that time, then we expect $x\in S$. If the trajectory does not reach the core, then we expect $x\not\in S$. Thus, by estimating the probability $\chi(x)$ to reach the core in a certain time-span $T$ starting in $x$ we get a membership function. This membership value is very close to a committor function value, if the core sets are assumed to be absorbing. The reason is, that a process starting in $S$ would quickly find the core. A process that needs ``too much'' time, probably found another core set and is trapped.  

Besides committor functions, there are other ideas to access $\chi$ efficiently. All methods which compute reaction coordinates and reaction paths \cite{reactionpaths} (as paths in high dimensional spaces) can also be used to approximate $\chi$, if we use the result of Sec.~\ref{sec:paths} that the holding probability decreases the most in the direction of that path. 
 
\subsection{Time discretization}
According to the ideas of the last section, one can find methods to  evaluate $\chi(x)$ which are not based on a linear combination of ansatz functions. The computation is based on reaction paths or on sampling from trajectories starting in $x$. Thus, the function evaluation $\chi(x)$ is a result of a simulation, i.e., of a discretization of time and not of a discretization of $\Gamma$ and, thus, circumvents the curse of dimensionality. We will apply this idea in order to find the paramters $\alpha$ and $\beta$ such that $\alpha\chi+\beta{1\!\! 1}$ approximates ${\cal L}^\ast \chi$ as good as possible.
However, even if it is possible to evaluate $\chi(x)$ pointwise by defining a core set and running simulations or, alternatively, by computing reaction paths starting in $x$, the pointwise computation of ${\cal L}^\ast \chi$ is not that straight forward. To solve this, we will exploit the fact that ${\cal P}^\tau=\exp(-\tau {\cal L}^\ast)$. Starting with the condition that we want to acchieve,
\[
  {\cal L}^\ast\chi=\alpha \chi +\beta,
\]
we get that ${\cal L}^\ast$ is a scale-shift-operator for $\chi$. For $i>0$, an iterative application of that operator leads to
\[
  \frac{\big(-\tau {\cal L}^\ast\big)^i}{i!}\chi=\frac{(-\tau\alpha)^i}{i!}\chi + \frac{(-\tau)^i\alpha^{i-1}\beta}{i!}.
\]
By taking the sum for $i=1,\ldots,\infty$ and adding $\chi$ on both sides
\begin{equation}\label{gammabasis}
  {\cal P}^{\tau}\chi=e^{-\tau\alpha}\chi + \frac{\beta}{\alpha}\big(e^{-\tau\alpha}-1\big).
\end{equation}

This has the following algorithmic consequences. Instead of ${\cal L}^\ast \chi$ we can evaluate ${\cal P}^\tau\chi$ pointwise: Given the point $x\in\Gamma$ for which we want to evaluate ${\cal P}^\tau\chi$, we start $M$ trajectories in $x$ of time-length $\tau$. For all the end points $x_\tau^{(k)}, k=1,\ldots,M$ of those trajectories we average over the values $\chi(x_\tau^{(k)})$. This provides the value of ${\cal P}^\tau\chi(x)$. We will also evaluate $\chi(x)$ at the starting point $x$ and, after that, solve the linear regression problem
\begin{equation}\label{regression}
   \min_{\gamma_1,\gamma_2}\|{\cal P}^\tau\chi(\cdot)-\gamma_1\chi(\cdot)-\gamma_2\|,
\end{equation} 
where every starting point $x$ generates one entry ${\cal P}^\tau\chi(x)-\gamma_1\chi(x)-\gamma_2$ of the vector. If the regression problem is exactly solvable, then  
\begin{equation}\label{gamma}
\gamma_1=e^{-\tau\alpha}, \quad \gamma_2=\frac{\beta}{\alpha}\big(e^{-\tau\alpha}-1\big)
\end{equation}
according to Eq. (\ref{gammabasis}). The proposed computation of the $\chi$-exit rate is depicted in Algorithm~\ref{algor}.
\begin{algorithm}
1. Determine a finite set $X$ of points in $\Gamma$.\vspace*{2mm}

2. For every point $x\in X$ evaluate $\chi(x)$ and ${\cal P}^\tau\chi(x)$ using simulations of (\ref{diffusion}).\vspace*{2mm}

3. Solve the linear regression problem (\ref{regression}) for the computation of $\gamma_1$ and $\gamma_2$.\vspace*{2mm}

4. With the aid of Eq. (\ref{gamma}) compute the values $$\alpha=-\frac{1}{\tau}\ln(\gamma_1)$$ and $$\beta=\frac{\alpha\gamma_2}{\gamma_1-1}.$$ 
5. The $\chi$-exit rate is $\epsilon_1=\alpha+\beta$.
\caption{\label{algor} Computing $\chi$-exit rates}
\end{algorithm}     

\subsection{Square-Root-Approximation of ${\cal L}^\ast$}
The sqrt-approximation is only needed for the artificial, illustrative examples in order to can compute analytical (non-statistical) results. A spatial discretization is, in principle, not needed for the application of the above theory. ${\cal L}^\ast$ is a continuous operator. For some easy examples, we will use a matrix  $L^\ast\in\RR^{n\times n}$ instead of ${\cal L}^\ast$. A possible, heuristic discretization scheme is available \cite{sqrt}. Note that $-L^\ast$ can be regarded as a transition rate matrix. If we assume a discretization of the state space $\Gamma$ into  $n$ subsets, then the transition rate between neighboring subsets $i$ and $j$ is given by $-L^\ast_{ij}=\sqrt{\pi_j/\pi_i}$, where $\pi_i=\exp(-\frac{1}{k_b T}V(i))$ is the Bolzmann weight and $V(i)$ is the potential energy value at the center of box $i$. $k_b$ is the Bolzman factor, $T$ is the temperature. The diagonal elements of $L^\ast$ are adjusted such that the row sum of $L^\ast$ is zero. This type of defining $L^\ast$ leads to a reversible process with a stationary distribution given by the Bolzmann distribution. We will use this square root approximation for the numerical examples below.

\section{Illustrative examples}
\subsection{Idea 1: $\chi$-exit rates from eigenvalues and eigenfunctions}\label{sec:idea1}
\begin{figure}[t]
\centering
\includegraphics[width=0.5\textwidth]{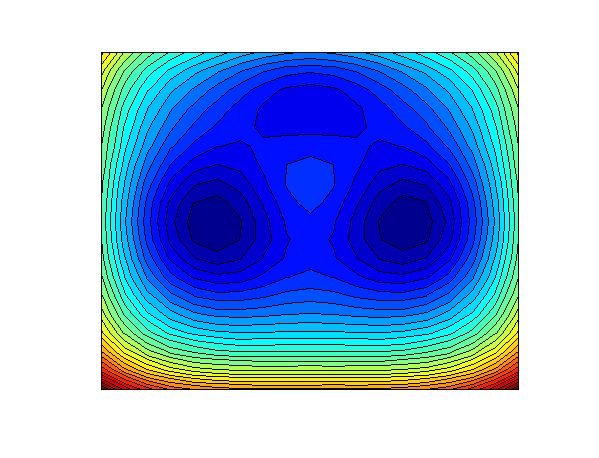}
\caption{\label{fig:potential} The potential energy function (\ref{eq:expo}) has two deep minima (dark blue, left and right) and a less deep minimum at the top. High values are indicated by red color, whereas low values are indicated by blue color.}
\end{figure}

The first example demonstrates how $\chi$-exit rates can be computed if a non-trivial eigenfunction and its eigenvalue $\bar{\epsilon}\not= 0$ of ${\cal L}^\ast$ are known. For this purpose the following $2$-dimensional potential energy function $V:\RR^2\rightarrow \RR$ is analyzed:
\begin{eqnarray}
\label{eq:expo}
V(x)&=&3\,\exp\big(-(4x_1-2)^2-(4x_2-{\frac{7}{3}})^2\big)\cr\cr &&-3\,\exp\big(-(4x_1-2)^2-(4x_2-{\frac{11}{3}})^2\big)\cr\cr
       &&-5\,\exp(-(4x_1-3)^2-(4x_2-2)^2)\cr\cr &&-5\,\exp(-(4x_1-1)^2-(4x_2-2)^2)\cr\cr
       &&+ 0.2\,(4x_1-2)^4+0.2\,(4x_2-{\frac{7}{3}})^4.
\end{eqnarray}
This function is depicted in Fig.~\ref{fig:potential}. In order to easily construct a discretized infinitesimal generator $-L^\ast$ on that potential, we generated a regular $50\times 50$ box discretization on the definition set $[0,1]\times [0,1]$. The transition rates between neighboring boxes $i$ and $j$ were set to be $-L^\ast_{ij}=\sqrt{\pi_j/\pi_i}$, where $\pi_i=\exp(-V(i))$ is the Bolzmann weight and $V(i)$ is the potential energy value at the center of box $i$. This construction is according to the proposed square root approximation of infinitesimal generators \cite{sqrt}.
\begin{figure}[ht]
\centering
\includegraphics[width=0.5\textwidth]{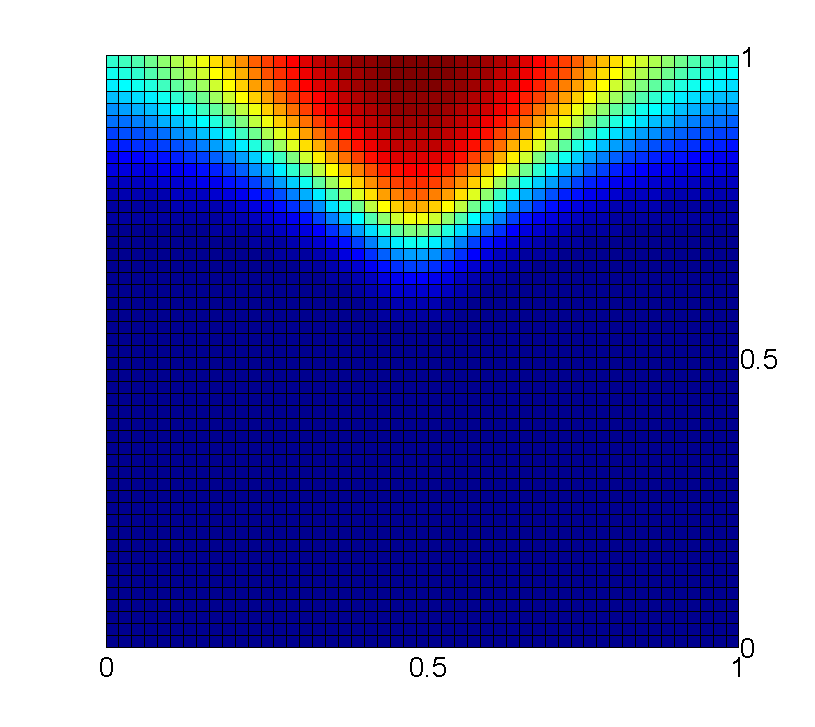}
\caption{\label{fig:eigf3} An approximation of one eigenfunction of ${\cal L}^\ast$ using the discretized operator $L^\ast$. This eigenfunction ``corresponds'' to the less deep minimum of the potential energy function.}
\end{figure}
With the aid of this discretization scheme, we can approximate one eigenfunction $f$, which has the third lowest eigenvalue $\bar{\epsilon}=0.0086$. This eigenfunction is shown in Fig.~\ref{fig:eigf3}. According to what has been said in Sec.~\ref{sec:paths} one can already derive the exit paths out of the top minimum of the potential energy surface from the negative gradient of the eigenfunction. Given a point $x$ in the definition set $\Gamma$, the holding probability decreases the most in the direction of $-\nabla f(x)$.  The highest value of the eigenfunction is $\max_i f(i)=0.0543$ whereas the minimal value is $\min_i f(i)=-0.0133$. 
Given one eigenfunction of ${\cal L}^\ast$, the computation of the membership function $\chi = \bar{\alpha} f +\bar{\beta}{1\!\! 1}$ based on PCCA+ is unique \cite{PCCAplus} with
\[
\bar{\alpha}=\frac{1}{\max_i f(i) - \min_i f(i)}
\]
and
\[
\bar{\beta}=\frac{-\min_i f(i)}{\max_i f(i) - \min_i f(i)}.
\]    
 
These quantities are sufficient to calculate the $\chi$-exit rate $\epsilon_1=\bar{\epsilon}(1-\bar{\beta})=0.0069$, also the statistical weight $\pi_\chi=\bar{\beta}=0.1965$, and the penalty parameter of the holding probability which is $\epsilon_2=\bar{\epsilon}\pi_\chi=0.0017$. The $\chi$-exit rate is physically meaningful according to Sec.~\ref{sec:meaning}, because $\epsilon_2\ll \epsilon_1$. 

\subsection{Comparison: fuzzy vs. set}
The mean holding time for a set $S$ is zero at the boundary of the set. The $\chi$-mean holding time $t_1(x)=\chi(x)/\epsilon_1$ is only zero for $\chi(x)=0$. Therefore, the set-based holding time $t(x)$ computed from the partial differential equation ${\cal L}^\ast t(x)=1$ must be different from $t_1(x)$. Note that in the molecular simulation setting, a simulated process has not reached the core of another conformation if the trajectory is at the boundary of the set $S$. If we want to compare the results of the $\chi$-mean holding time with a set-based approach, the fuzzyness of the answer to the question ``where in $\Gamma$ do we reach another conformation?'' plays an important role. If we, e.g., define $S$ in the situation of Sec.~\ref{sec:idea1} to be that part of $\Gamma$ which is defined by 
\[
S=\{x\in\Gamma \,|\, \chi(x)>0.22\},
\]  
then a set-based holding probability will be zero at the boundary of $S$, while the $\chi$-holding probability will be $\frac{0.22}{\epsilon_1}=31.88$. This is a huge difference.
\begin{figure}[ht]
\centering
\includegraphics[width=0.5\textwidth]{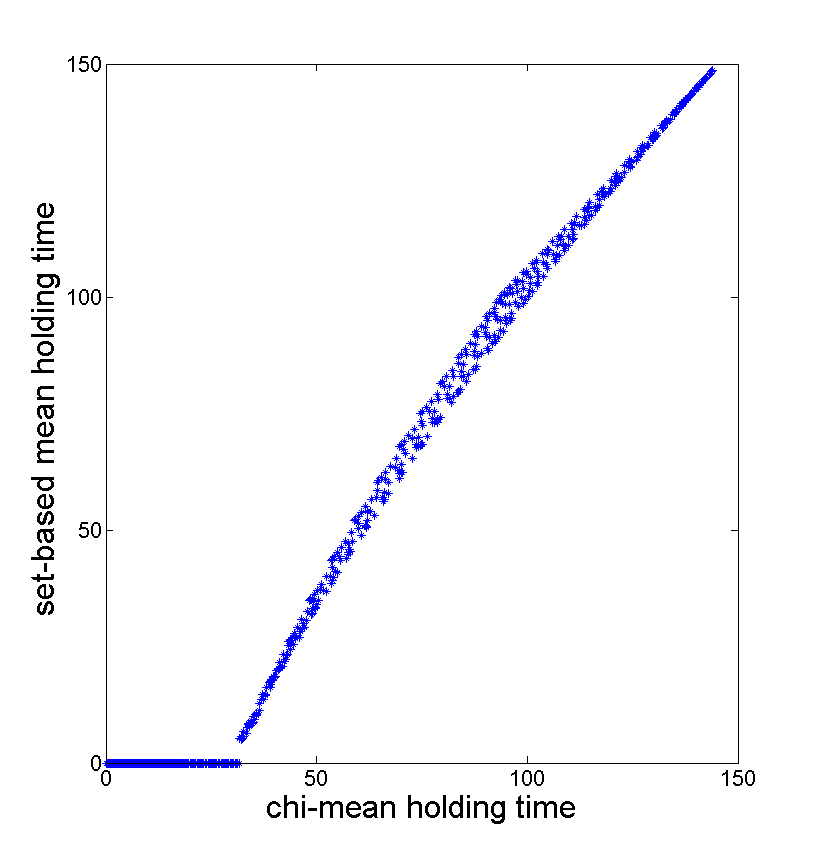}
\caption{\label{fig:comparison} In this plot the $x$-axis is the $\chi$-mean holding time $t_1(x)$ and the $y$-axis is the set-based mean holding time $t(x)$ compared for the values of the $2500$ discretization boxes in the situation of Sec.~\ref{sec:idea1} and for  $S=\{x\in\Gamma \,|\, \chi(x)>0.22\}$. The two quantities are mainly identical, except for a different behavior at the boundary of $S$.} 
\end{figure}
In Fig.~\ref{fig:comparison}, we compare the $\chi$-mean holding time with the computed set-based mean holding time for the $2500$ cells of the discretization of ${\cal L}^\ast$. There is a clear correlation between these two quantities, except for the fact, that the $\chi$-mean holding time $t_1(x)$ does not have the zero-plateau. For a definition of $S$ on the basis of the condition $\chi(x)>0.5$, the set-based holding time $t(x)$ would be much smaller than $t_1(x)$. For a definition of $S$ on the basis of the condition $\chi(x)>0.1$, the set-based holding time would mostly be much higher than $t_1(x)$. The $\chi$-mean holding time is, thus, like a ``compromise'' in that sense. It is like a ``mean'' mean holding time for different possible choices of the boundary of the metastable set $S$.     

\subsection{Idea 2: Linear regression of PCCA+ results}\label{sec:idea2}

\begin{figure}[ht]
\centering
\includegraphics[width=0.5\textwidth]{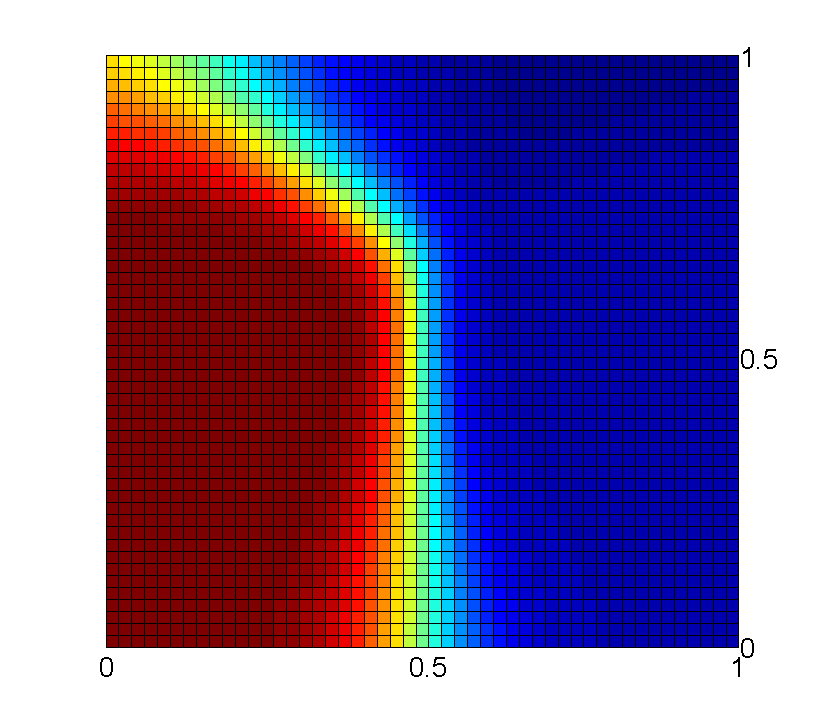}
\caption{\label{fig:chi2} The membership function $\chi$ which represents the starting point assignment for the left deep minimum of the potential energy surface of Fig.~\ref{fig:potential}. This function is a linear combination of three eigenfunctions of ${\cal L}^\ast$.}
\end{figure}

The case, that has been described in Sec.~\ref{sec:idea1} is artificial. Usually, not every single eigenfunction of ${\cal L}^\ast$ can be interpreted as a membership function. Sometimes the membership function $\chi$ has to be composed as the linear combination of several eigenfunctions. This situation is shown in Fig.~\ref{fig:chi2}.
In this situation 
\[
\chi= 0.4452 \cdot{1\!\! 1} + 17.7865\cdot f_2 - 4.1266\cdot f_3,
\]
where $f_2$ and $f_3$ are approximated eigenfunctions of ${\cal L}^\ast$ corresponding to the eigenvalues $0.0025$ and $0.0086$, respectively. From the first factor, we can directly extract $\pi_\chi=0.4452$ which is less than $0.5$, thus, we will have the good case that $\epsilon_2 < \epsilon_1$. In the given situation it is possible to compute ${\cal L}^\ast\chi$ analytically: 
\[
{\cal L}^\ast\chi = 17.7865 \cdot 0.0025\cdot f_2 - 4.1266 \cdot 0.0086\cdot f_3.
\]    
In the case of several eigenfunctions, ${\cal L}^\ast\chi$ is not a linear combination of ${1\!\! 1}$ and $\chi$ any more. But if we solve the linear regression problem of minimizing the norm $\|{\cal L}^\ast\chi -\alpha\chi - \beta {1\!\! 1}\|$, then the result is $\alpha=0.0028$ and $\beta= -0.0014$. Thus, the $\chi$-exit rate is $\epsilon_1=\alpha+\beta=0.0014$.

\subsection{Idea 3: Committor functions as approximation space}\label{sec:idea3}

\begin{figure}[t]
\centering
\includegraphics[width=0.5\textwidth]{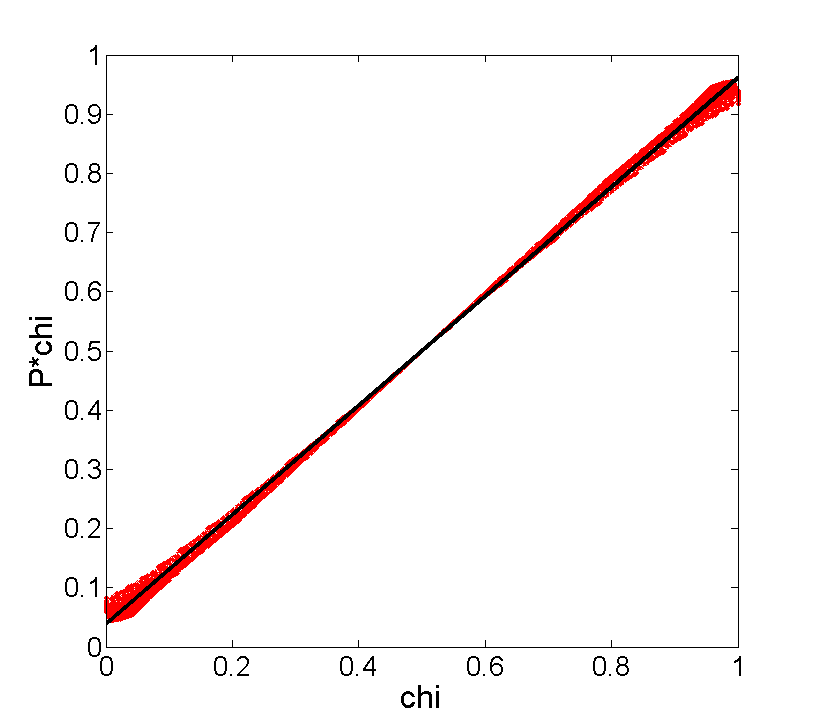}
\caption{\label{fig:chi_Pchi_exact} Computing $\chi$ as a committor function of the process and $P^\tau\chi$ for $\tau=100$. The red points correspond to the plot of the $2500$ $\chi$-values against their propagated values $P^\tau\chi$. The black line is the $\|\cdot\|_2$- regression result. Many possible lines can fit the given data depending on the regression norm.}
\end{figure}

Given the square root approximation $L^\ast$ (with $k_b T=1$) one can easily compute the discretized committor function of the process between the left and the right deep minimum. This committor function will now serve as an approximation for the membership function, it will, therefore, also be denoted as $\chi$. For the committor function between the two minima, two core sets are needed. The core sets are based on the $50 \times 50$ discretization of $\Gamma$. All discretization boxes having a statistical weight higher than $0.0025$ are assigned to one of the core sets. After computing $\chi$, the propagated values $P^\tau\chi$ with $\tau=100$ are determined. Note that $P^\tau=\exp(-\tau L^\ast)$. In Fig.~\ref{fig:chi_Pchi_exact} the $\chi$-values are plotted against the $P^\tau \chi$-values for solving the regression problem (\ref{regression}). Depending on the regression norm, the regression results may differ. Taking the $\|\cdot\|_2$-norm, the results are $\gamma_1=0.8201$ and $\gamma_2=0.900$. Thus, $\alpha=0.0020$ and $\beta =  -0.0010$. The $\chi$-exit rate is approximated to be $\epsilon_1=0.0010$, which is lower than the result of Sec.~\ref{sec:idea2}.        

\subsection{Idea 4: Short-time simulations for estimating committor functions}
\label{sec:idea4}

\begin{figure}[ht]
\centering
\includegraphics[width=0.5\textwidth]{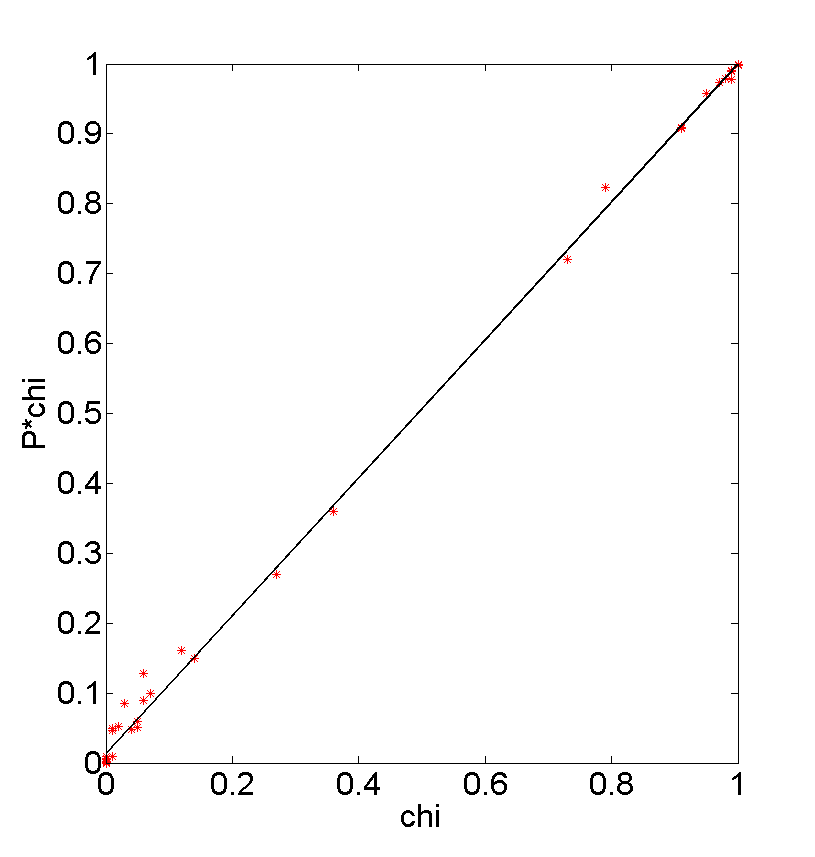}
\caption{\label{fig:Pchi} ``Evaluating'' $\chi$ and ${\cal P}^\tau\chi$ at $50$ randomly chosen points $x\in \Gamma$. The relation betwen these two quantities should be linear. By linear regression a line is fitted into the data providing $\gamma_1=0.9870$ and $\gamma_2=0.0128$.}
\end{figure}
In the three illustrative examples above, a discretized version of the infinitesimal generator was given. Thus, we were able to compute (approximative) eigenfunctions and eigenvalues of ${\cal L}^\ast$ that can be used to define membership functions and $\chi$-exit rates out of those corresponding fuzzy sets $\chi$. The need for a discretization is a drawback of that method in high-dimensional conformational spaces of molecular systems. In Sec.~\ref{sec:complex} it has been discussed, whether there is a way to estimate the (committor) function $\chi$ by running simulations only. If this is possible, then the Algorithm~\ref{algor} would also provide $\chi$-exit rates. In order to illustrate how this algorithm practically works, we took again the potential energy function of Fig.~\ref{fig:potential}. This time we will apply a Brownian dynamics simulation of ${\cal P}^\tau$ instead of the square root approximation of ${\cal L}^\ast$. The Brownian dynamics is chosen to be ``faster'' than the square-root approximation: For computing the value of $\chi(x)$ for each $x\in X$, we start $100$ trajectories in $x$ according to (\ref{diffusion}) with $\sigma = 0.8$ and an Euler-Maruyama time discretization of $\delta t=0.001$. $\chi(x)$ is defined as the percentage of those trajectories which have reached a certain core set within less than $100$ integration steps. The core set is reached, if the $x_1$-coordinate is in the interval $[0.2, 0.3]$ and the $x_2$-coordinate is in the interval $[0.4, 0.5]$. By this procedure we get a $\chi$-function which is very similar to that in Fig.~\ref{fig:chi2}. Different from the situation of Sec.~\ref{sec:idea2}, $\chi$ is only given point-wise. $50$ randomly chosen points $x$ in the box $[0,1]^2$ have been used for this $\chi$-function evaluation. ${\cal P}^\tau\chi(x)$ has also been evaluated in a similar way. $100$ trajectories with $50$ integration steps only (total time length is $\tau = 0.05$) have been generated to propagate $x$. At those propagated points $x^{(k)}$ the $\chi$-function has been evaluated as described above. ${\cal P}^\tau\chi(x)$ is given by the averaged $\chi$-value at the propagated points. In Fig.~\ref{fig:Pchi} the ${\cal P}^\tau\chi$-values are plotted against the $\chi$-values at the $50$ chosen points. 

If $\chi$ had been a linear combination of an eigenfunction of ${\cal L}^\ast$ and ${1 \!\! 1}$, then this plot would show a line. From the axis intercept $\gamma_2$ and the slope $\gamma_1$ of that line, the $\chi$-exit rate is estimated. Small deviations from the perfect line due to sampling errors, however, lead to high relative errors in $\epsilon_1$. By fitting a line to the data points in Fig.~\ref{fig:Pchi}, we estimated the $\chi$-exit rate $\epsilon_1=0.0042$.

\section{Molecular example}\label{sec:molecule}

\begin{figure}[hb]
	\centering
		\begin{tikzpicture}[rotor/.style={inner sep=3pt,outer sep=0,fill opacity=1,minimum width=2cm,circle}]
		%rotors
		\node[rotor] (n1) at (0,0.2) {};
		\node[rotor] (n2) at (1.1,-0.25) {};	
		\node[rotor] at (-1.1,0) {$\phi$};
		\node[rotor] at (0.3,0.8)  {$\psi$};
		\node(c1) at (-3,0) {C};
		\node(c2) at (-1.5,0.8)  {C};
		\node(c3) at (0,0) {C};
		\node(c4) at (1.5,0.8)  {C};
		\node(c5) at (3,0) {C};
		\draw[line width=1pt] (c1)--(c2)--(c3)--(c4)--(c5);
		% rotation direction arrows
		\draw [<-,line width=1pt] (n1) ++(140:5mm) --++(-30:-1pt) arc (50:270:3mm);
		\draw [<-,line width=1pt] (n2) ++(140:5mm) --++(200:-1pt) arc (-70:150:3mm);
		\end{tikzpicture}
	\caption{\label{angles}The pentante molecule consists of $5$ carbon atoms and $12$ hydrogen atoms (not shown).  Two torsional angles $\phi$ and $\psi$ determine the conformation of that molecule. The depicted conformation corresponds to $\phi=\psi=180^\circ$, which is the most stable conformation. For this conformation we aim to compute the exit rate.}
\end{figure}
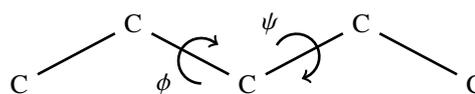

Algorithm~\ref{algor} can easily be applied to molecular systems. We will demonstrate this algorithm for the simulation of an n-pentane molecule, shown in Fig.~\ref{angles}. This molecule has $17$ atoms.  Thus, its state space $\Gamma=\RR^{3\cdot 17}$ is $51$-dimensional.\\ 

\begin{figure}[htb]
	\centering
    \includegraphics[width=0.5\textwidth]{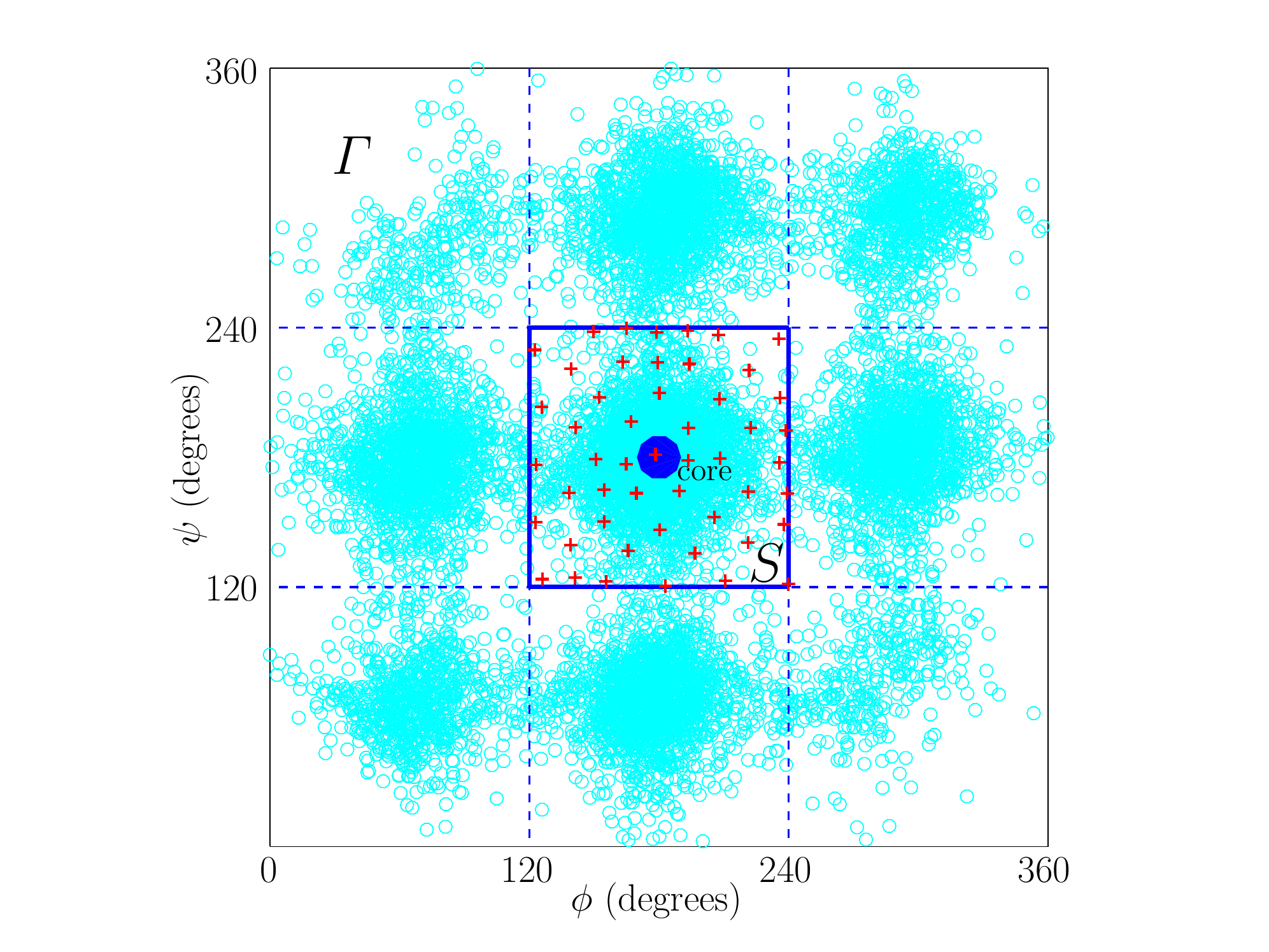}
	\caption{\label{core} Torsion angle distribution of the two torsion angles of pentane at 700K, for the interval of 0 to 360 degrees (cyan circles). Sampling starting points $x\in X$ for the Algorithm~\ref{algor} from the box $[120,240]^2=:S$ (red crosses). Total number of starting points is $50$. } 
\end{figure}

{\bf Defining the core set.}
The definition of the core set depends on the kind of rare event statistics which we want to estimate. If we want to figure out the slow diffusion of the pentane molecule in the $3$-dimensional space, then we would project the $51$-dimensional coordinates, e.g., onto the center of mass of pentane and define the core set as a ball in this $3$-dimensional space. However, chemists are more interested in the {\em internal} transitions of the molecule (i.e., its conformations). It is well-known that the pentane molecule has several metastable low-energy conformations. Those can be determined by considering two torsional angles $\phi$ and $\psi$. Each torsion angle is defined by $4$ consecutive carbon atoms, see Fig~\ref{angles}. A long-term simulation of pentane at $700$K (with fast transitions between the conformations) reveals that there are nine different peaks of the stationary (Boltzmann) distribution in the $\phi$-$\psi$-digram shown in Fig.~\ref{core} (cyan circles and dotted boxes).
If we want to estimate the exit rate from the central conformation of pentane, then the core set can be defined as all states $x\in\Gamma$ of pentane which have a pair of torsion angles approximately at $(\phi,\psi)\approx (180,180)$. This core set is indeed a non-convex, connected, unbound set in $\Gamma=\RR^{51}$. In the $\phi$-$\psi$-plane it is a circle, see the blue circle in the center of Fig.~\ref{core}. Although the described core sets can easily be  projected to a $3$- or $2$-dimensional space, neither the potential energy $V$ nor the membership function $\chi$ is a $3$- or $2$-dimensional function. This example of pentane is indeed $51$-dimensional with $\chi:\RR^{51}\rightarrow [0,1]$.\\     

\begin{figure}[htb]
	\centering
	\includegraphics[width=0.5\textwidth]{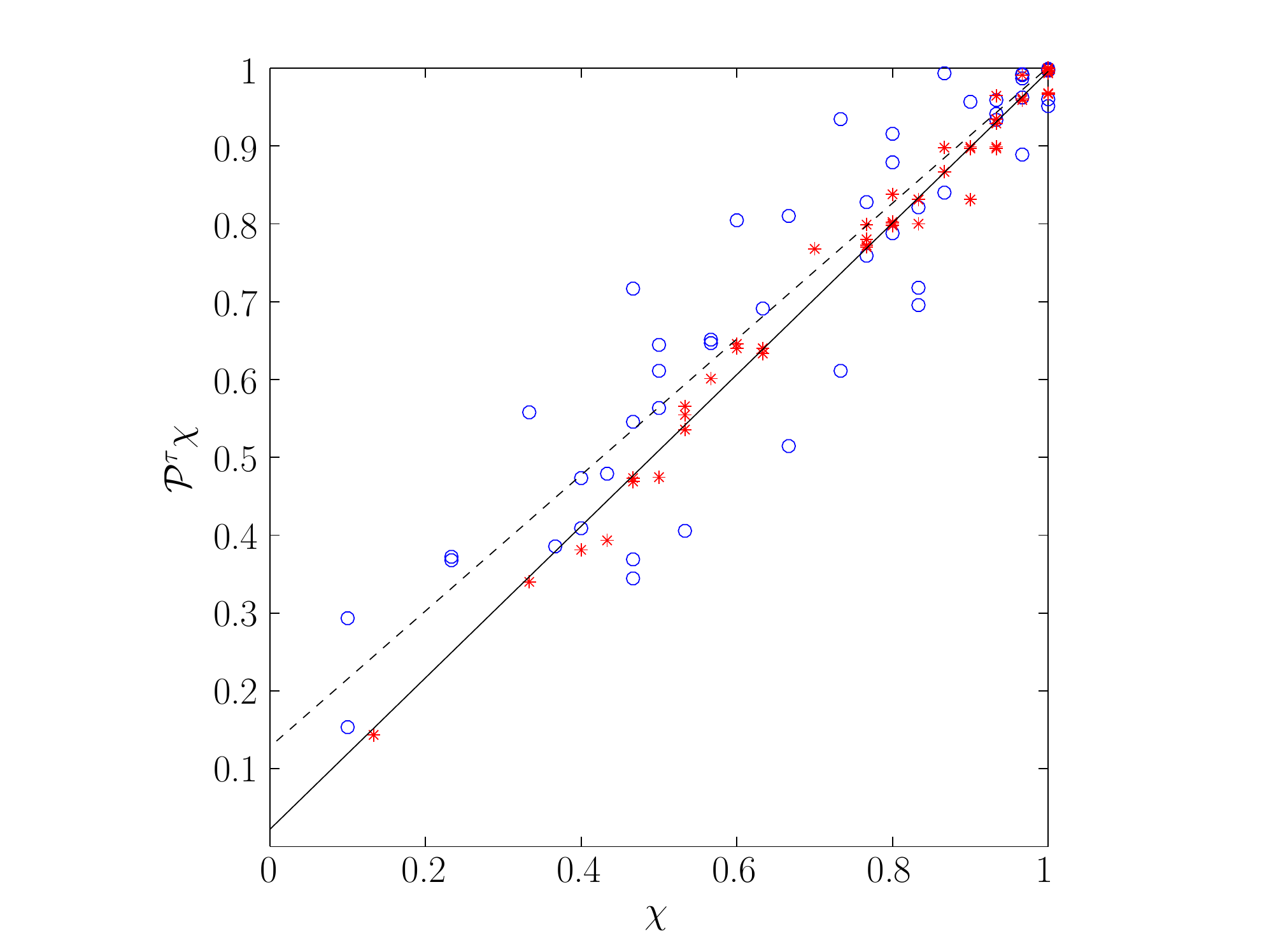}
	\caption{\label{lin_regr} $\chi$ and $\mathcal{P}^\tau\chi$ for the indicated $50$ starting points in Fig.~\ref{core}. The red stars correspond to a GROMACS SD simulation. The linear regression of the corresponding $\chi$-${\cal P}^\tau \chi$-correlation is given by a solid line. If an MD simulation is performed (blue circles), then the correlation between $\chi$ and ${\cal P}^\tau\chi$ is worse (dotted line) and does not lead to an interpretable rate estimation.}
\end{figure}

{\bf Application of Algorithm~\ref{algor}.}
Instead of discretizing the $51$-dimensional space for approximating $\chi$, which would lead to a curse of dimensionalty, we will only evaluate $\chi(x)$ and ${\cal P}^\tau \chi(x)$ at $50$ different points $x\in X$. The projection of those $50$ points onto the $\phi$-$\psi$-plane is shown as red crosses in Fig.~\ref{core}. 
To compute the $\chi(x)$-value 30 GROMACS stochastic dynamics (SD) simulations for each $x\in X$ with time discretization step $\delta t = 0.001$ps at a temperature of $310$K (conformational transitions are rare events at this temperature) were performed.  All molecular simulations are performed with GROMACS 5.1.2 \cite{GROM, Gro1}. Based on this simulation the percentage of the trajectories which have reached the core within first 0.5ps was calculated. The coordinates at the end points of these simulations are used for computing the value of $\mathcal{P}^\tau \chi$. For each point 30 more simulations with 1000 time steps and the same time step $\delta t = 0.001$ps (total time is 1ps) were done. $\mathcal{P}^\tau \chi (x)$ was calculated as the average of the $\chi$-values for the propagated points. Using linear regression we get $\gamma_1=0.9738$, $\gamma_2=0.0220$. Therefore, the $\chi$-exit rate is computed to be $\epsilon_1=0.01ps^{-1}$. This result means, that on average it needs about $100$ps to exit from the central conformation of n-pentane. However, to yield this result, we used only much shorter trajectories (of $0.5$ps or $1$ps) which could easily be generated independently in parallel on different processors. \\

{\bf MD versus SD.}
For showing that the theory depends on a certain kind of equations of motion, we performed molecular dynamics (MD) simulations in which all parameters and starting states were chosen to be equal to the SD simulation. MD simulation is based on Newton’s equations of motion and not on (\ref{diffusion}). For detailed information, please, see chapter 3.8 of the GROMACS manual \cite{GROM, Gro1}. In the MD case, the linear regression led to a different result with $\gamma_1=0.8731$ and $\gamma_2=0.1280$ (see Fig.~\ref{lin_regr}). Thus, $\alpha=0.2714$ and  $\beta=-0.2738$. Therefore, we get a negative value for the exit rate. For the presented theory it is mandatory to use (\ref{diffusion}) as the equations of motion. MD is not applicable.\\

{\bf Parallelization.}
 For our approach, two different kinds of computational parallelization can be combined. There are already some built-in parallelization schemes in GROMACS to run trajectories. As we additionally know that all 50 start states are independent, it is trivial to parallelize the simulation of those trajectories as well. For the MD/SD simulations and calculation of $\chi$, $\mathcal{P}^\tau \chi$ values all available cores were used. In order to achieve optimal simulation performance the number of starting trajectories should be divisible by the number of used cores.\\

\begin{figure}[htb]
	\centering 
	\includegraphics[width=0.5\textwidth]{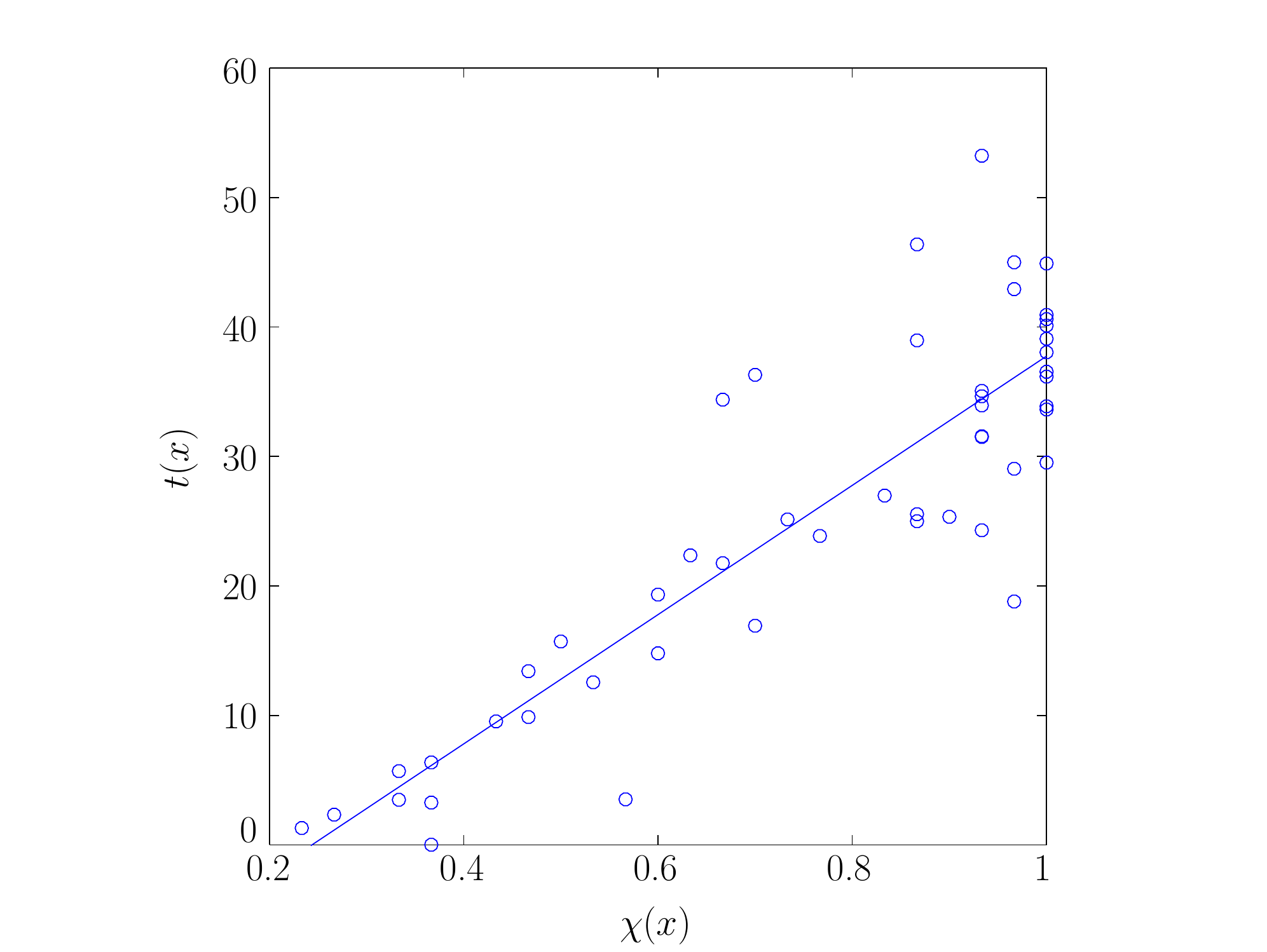}
	\caption{\label{chi_tau} $\chi$ and $\tau$ for the indicated $50$ starting points in Fig.~\ref{core}. There is a relation between the membership value $\chi(x)$ and the simulated mean exit time $t(x)$.}
\end{figure}

{\bf Validation.}
For a validation we will compare the $\chi$-exit rate with a set-based exit rate. For this validation we, thus, need to define a starting set $S$, which makes this comparision difficult, because the quality of our result will depend on the choice of $S$. We computed the mean first exit time $t$ for leaving the blue box $S$ indicated in Fig.~\ref{core}. This is an {\em arbitrarily} chosen starting set according to chemical intuition about the dominant conformation of $n$-pentane. To estimate the exit rate, we again started $30$ SD-simulations from the $50$ starting points given as red crosses in Fig.~\ref{core}. These simulations were $200$ times longer ($100$ps) than our simulations for the $\chi$-exit rate estimation.  In principle, by this estimation a function $t(x)$ of set-based mean holding times in $\Gamma$ is approximated point-wise.  From Fig.~\ref{chi_tau} we see that the starting points with $\chi(x)=1$ have a mean holding time of about $40$ps. This would mean, that the set-based exit rate should be about $0.025$ps$^{-1}$, which is in the same order of magnitude as our result $0.01$ps$^{-1}$ but $2.5$ times higher. According to our presented theory, we additionally know that $t_1(x)=\frac{1}{\epsilon_1}\chi(x)$ is the $\chi$-mean holding time. Like in Fig.~\ref{fig:comparison}: If $t_1(x)$ would be a good representation of the set-based mean holding time behavior $t(x)$, then there should be a linear dependence between the simulated exit time and the membership value $\chi(x)$ as it is indeed the case in Fig.~\ref{chi_tau}. Note that the $t(x)$-computation suffers from high variance for larger holding times.  

We calculated the mean exit time for each starting point based on the 30 SD-simulations. Only 1 of 50 starting points has the mean exit time less than 1ps, therefore we can say that in this case the exit rate is equal to $\tau=0.02$.

%%%%%%%%%%%%%%%%%%%%%%%%%%%%
\section{Conclusion}  
If one question is too complicated to answer, then maybe there is a slightly different question which provides the same kind of technical value and reveals simple relations.  We started the article with the observation that a potential-driven diffusion process tends to sample from the core of a metastable set $S$. Trajectories starting in some point $x\in S$ which are constructed according to Eq.~(\ref{diffusion}) probably quickly reach the core and extremely rarely leave the set $S$. In Sec.~\ref{sec:complex} we inverted this observation and turned it into a definition of a membership function $\chi$. If a trajectory starts in $x\in \Gamma$ and quickly reaches a pre-defined core, then we expect that $x$ is part of the (implicitly pre-defined) metastability. The function value $\chi(x)$ is defined as the portion of trajectories which start at $x$ and reach the core quickly. A function $\chi$ of that type can be efficiently estimated point-wise by running simulations. There is no curse of dimensionality, if the $\chi$-function is approximated point-wise. Given that the metastability is now a fuzzy set $\chi$, then we need  a new definition of what we want to understand to be the {\em holding probability}. In this article we introduced a definition of a $\chi$-holding probability which provides easy mathematical relations between $\chi$-exit rates, $\chi$-exit paths, and $\chi$-mean first exit times.   In principle, the $\chi$-exit rate is given by the slope and the axis intercept of the linear relation between $\chi$ and ${\cal P}^\tau \chi$. For recovering the linear relation between $\chi$ and ${\cal P}^\tau \chi$, it is very ``useful'' that $\chi$ has values ranging from $0$ to $1$, i.e., the information of ``how much does a state $x$ belong to the metastability'' is exploited.

{\bf Acknowledgement.} The work has partially been financed by the CRC-1114 ``Scaling Cascades in Complex Systems'', project A05. Data and m-scripts can be found at \url{http://www.zib.de/ext-data/soft_exit_rate}.

\end{document}